\documentclass[11pt, twoside]{amsart}

\usepackage[utf8]{inputenc}
\usepackage[T1]{fontenc}
\usepackage{amssymb,amsmath,amstext}
\usepackage{hyperref}

\newtheorem{theor}{Theorem}[section]
\newtheorem{lem}[theor]{Lemma}
\newtheorem{defin}[theor]{Definition}

\newtheorem{prop}[theor]{Proposition}

\newtheorem{cor}[theor]{Corollary}
\newtheorem{rem}[theor]{Remark}

\newtheorem{fact}[theor]{Fact}

\numberwithin{equation}{section}

\newcommand{\acl}{\mathrm{acl}}
\newcommand{\dcl}{\mathrm{dcl}}
\newcommand{\tp}{\mathrm{tp}}

\newcommand{\cl}{\mathrm{cl}}

\newcommand{\es}{\emptyset}

\newcommand{\su}{\mathrm{SU}}

\newcommand{\nts}{\negthickspace}
\newcommand{\uhrc}{\nts \upharpoonright \nts}

\newcommand{\meq}{^{\mathrm{eq}}}

\newcommand{\mcA}{\mathcal{A}}

\newcommand{\mcM}{\mathcal{M}}
\newcommand{\mcN}{\mathcal{N}}

\newcommand{\mrqf}{\mathrm{qf}}

\newcommand{\ind}{\raisebox{-2pt}[5pt][0pt]{$\smile$} \hspace*{-6.8pt}\raisebox{3pt}[5pt][0pt]{$|$} \; \: }
\newcommand{\nind}{\raisebox{-2pt}[5pt][0pt]{$\smile$} 
\hspace*{-6.8pt}\raisebox{3pt}[5pt][0pt]{$|$}\hspace*{-6.8pt}
\raisebox{3pt}[5pt][0pt]{$\diagup$} }

\newcommand{\rng}{\mathrm{rng}}

\title[Supersimple $\omega$-categorical theories]
{Supersimple $\omega$-categorical theories and pregeometries}

\author{Vera Koponen}

\address{Vera Koponen, Department of Mathematics, Uppsala University, Box 480,
75106 Uppsala, Sweden.}

\email{vera.koponen@math.uu.se}

\date{29 June 2018, revised version}

\begin{document}

\begin{abstract}{
We prove that if $T$ is an $\omega$-categorical supersimple theory with nontrivial dependence (given by forking), 
then there is a nontrivial regular 1-type over a finite
set of reals which is realized by real elements; hence forking induces a nontrivial pregeometry on the solution set of this type
and the pregeometry is definable (using only finitely many parameters).
The assumption about $\omega$-categoricity is necessary.
This result is used to prove the following: If $V$ is a finite relational vocabulary 
with maximal arity 3 and $T$ is a supersimple $V$-theory with elimination of quantifiers, 
then $T$ has trivial dependence and finite SU-rank.
This immediately gives the following strengthening of~\cite[Theorem~4.1]{Kop17a}:
if $\mcM$ is a ternary simple homogeneous structure with only finitely many constraints, then $Th(\mcM)$ has trivial dependence and
finite SU-rank.

}
\end{abstract}

\maketitle

\section{Introduction}

\noindent
The idea that ``global'' properties of a theory can, under some circumstances, be understood to a large part by its 
``local'' properties dates back, at least, to Zilber's studies of uncountably categorical theories in the 80'ies.
(For a monograph in English on this topic see \cite{Zilber_mono}.)
Ever since, this idea has been an important guideline in model theory giving rise to many
results, in particular in the studies of stable theories and later and
more generally in the studies of simple theories.

By a local property we mean a property of (the set of elements realizing) a type (i.e. a consistent set of formulas).
Certain types are of particular interest. 
Zilber considered minimal types and Hrushovski \cite{Hru85, Hru90, Hru92} 
generalized many of the results to regular types, in both cases in the context of stable theories.
Later such results where generalized to the context of simple theories.
If $T$ is a simple theory and $p$ is a regular type, then dividing (or equivalently forking) dependence, induces 
a pregeometry (or matroid) on the set of elements that realize $p$.
The interesting thing about this is that
the complexity of the pregeometries on regular types tend to reflect the complexity of dividing dependence of the theory.
This is, for example, the case when if one considers supersimple 1-based structures with finite SU-rank
as shown by Hart, Kim and Pillay~\cite{HKP} and by de Piro and Kim~\cite{DK}.
What is more directly relevant for this article is the following result of Goode~\cite[Propositions~2 and~5]{Goode}:
if $T$ is superstable with finite SU-rank (usually called U-rank in the context of stable theories) and
$T$ has nontrivial dependence (Definition~\ref{definition of trivial dependence}), then there is a nontrivial regular type
(Definition~\ref{definition of minimal dependent set}).
Palac\'{i}n~\cite{Pal17} observed that Goode's proof generalizes to supersimple theories with finite SU-rank.
On the other hand the assumption that the SU-rank is finite is necessary (even in the stable context) 
which is demonstrated by an example in~\cite{Goode}
(see Remark~\ref{remark about the use of omega-categoricity} below).

However, if we assume that the theory is $\omega$-categorical, then the assumption about finite SU-rank can be removed, which
is our first main result (Theorem~\ref{nontrivial dependence implies a nontrivial regular type}). More precisely: 
if $T$ is supersimple and $\omega$-categorical with nontrivial dependence, 
then there are $\mcM \models T$, a finite set $C \subseteq M$ and a nontrivial regular 1-type $p$ over $C$ realized by real elements.
In our application of this result it matters that $p$ is realized by real elements, as opposed to properly imaginary ones.
A result like this can be useful to show that certain theories cannot have complicated behaviour of dependence,
by showing that they cannot accomodate a complicated definable pregeometry.
In this article we consider theories with elimination of quantifiers in a finite relational language. 
Such theories are $\omega$-categorical by a well known characterization of $\omega$-categorical theories.
To facilitate this approach we prove the following  (Theorem~\ref{nontrivial dependence implies nontrivial pregeometry}):
If $T$ is a supersimple $V$-theory with elimination of quantifiers, where $V$ is a finite and relational vocabulary and $T$ has nontrivial dependence,
then there is another finite relational vocabulary $V'$ of the same maximal arity as $V$ and a 
supersimple $V'$-theory $T'$ with elimination of quantifiers such that dependence induces a nontrivial pregeometry on any model of $T'$
and $T'$ has a unique 1-type over the algebraic closure of $\es$ with imaginaries.
Then we show that no such $T'$ can exist if the maximal arity of $V'$ is 3 (Theorem~\ref{impossibility of a structure satisfying the other result}).
Hence, every supersimple $T$ with elimination of quantifiers in a finite relational language with maximal arity 3 (of the relation symbols)
has trivial dependence. 
It follows from work of Macpherson~\cite{Mac91} and de Piro and Kim~\cite{DK} 
that if $T$ is simple with elimination of quantifiers in a finite relational langauge
and $T$ has {\em non}trivial dependence, then $T$ is {\em not} 1-based, so if $T$ is supersimple then it must have a
{\em non}modular regular type.
This seems very difficult to achieve (with elimination of quantifiers in a finite relational langauge) so my guess is that 
dependence is always trivial in such a theory.
As a contrast, if we forget about elimination of quantifiers,
Hrushovski has constructued, with only a ternary relation symbol, a supersimple $\omega$-categorical structure
with SU-rank~1 and nontrivial (even not 1-based) independence 
(see for instance \cite[Chapter~6.2]{Wag} or \cite[Chapter~6]{Kim_book}).

Besides the nature of regular types and pregeometries, the SU-rank gives important information about a supersimple theory.
Every $\omega$-categorical superstable theory is $\omega$-stable as proved  by Lachlan~\cite{Lach74} and has finite Morley-rank
as proved by Cherlin, Harrington and Lachlan~\cite{CHL}, 
from which it follows
that it has finite SU-rank. 
There is a conjecture \cite[p.~205]{Wag} that every $\omega$-categorical supersimple theory has finite SU-rank, but the cases in
which it has been proved are still limited. 
For every smoothly approximable structure its theory is $\omega$-categorical and supersimple with finite SU-rank~\cite{CH, HKP}.
More generally, Evans and Wagner proved that every theory which is $\omega$-categorical, supersimple and CM-trivial has finite SU-rank
\cite{EW}.
If $T$ is simple with elimination of quantifiers in a finite relational language with maximal arity 2,
then $T$ is supersimple with finite SU-rank~\cite{Kop16PAM}.
We show here, in Corollary~\ref{trivial dependence and finite rank},
that if $T$ is supersimple with elimination of quantifiers in a finite relational language with maximal arity 3,
then $T$ has finite SU-rank.
This is actually a consequence of the mentioned result that such a theory has trivial dependence and the following result,
due to Palac\'{i}n~\cite{Pal17}:
If $T$ is $\omega$-categorical and supersimple with trivial dependence, then $T$ has finite SU-rank.
From Corollary~\ref{trivial dependence and finite rank}
we can directly improve the main result in~\cite{Kop17a} to the following:
If $\mcM$ is ternary, homogeneous (definition follows below) and simple with only finitely many constraints, 
then its theory is supersimple with finite SU-rank and trivial dependence.
To give some more background, Lachlan~\cite{Lach74} once
conjectured that every $\omega$-categorical stable theory is superstable, but this was disproved by 
Hrushovski~\cite{Hru_pseudoplane}.\footnote{
His counterexample of an $\omega$-categorical stable ``pseudoplane'' which is not superstable
has never been published as far as I know, but it is mentioned in various places in the literature.}

In some sense the simplest $\omega$-categorical theories are those with elimination of quantifiers, and yet, via 
so-called amalgamation constructions, it is possible to construct uncountably many supersimple theories with SU-rank~1 and with
elimination of quantifiers by using only one ternary relation symbol~\cite{Kop17a}.
Then it is natural to first try to solve hard problems in this special case, which is still challenging.
Moreover, theories with elimination of quantifiers, and in particular their countable models, are interesting from other perspectives.
Suppose that $T$ is a complete $V$-theory where $V$ is a finite relational vocabulary and let $\mcM$ be a {\em countable} model of $T$.
Then the following three conditions are equivalent 
(see for example~\cite[Chapter~7]{Hod}): (a) $T$ has elimination of quantifiers; (b) $\mcM$ is {\em homogeneous}\footnote{
The terminology {\em ultrahomogeneous} and {\em finitely homogeneous} is also used in the literature.}, 
meaning that every isomorphism between finite substructures of $\mcM$ can be extended to an automorphism of $\mcM$; 
(c) the class of finite $V$-structures that can be embedded into $\mcM$ has the amalgamation property and consequently
$\mcM$ is the so called Fra\"{i}ss\'{e} limit of this class. 

Homogeneous structures have turned out to be interesting from a variety of viewpoints. Since they have a rich automorphism
group they are interesting in permutation group theory and there are, via structural Ramsey theory, interesting connections
to topological dynamics.
Homogeneous structures have also become an important object of study in the area of constraint satisfaction problems.
See the survey article~\cite{Mac11} by Macpherson for more about the mentioned aspects of homogeneous structures and further references.
All stable homogeneous structures are well understood through the work of Lachlan and others~\cite{Lach97},
but new challenges are offered in the broader class of simple homogeneous structures. (For example, every stable
homogeneous structure is finitely constrained, but given only a ternary relation symbol one can construct uncountably many
{\em not} finitely constrained simple homogeneous structures~\cite{Kop17a}.)

The idea to understand (super)simple theories with elimination of quantifiers by considering a bound on the maximal arity of the
relation symbols may seem futile as the bound can always be increased.
But I suspect that once the case of maximal arity 4 is understood with respect to general questions such as triviality of
dependence and SU-rank, then we also get the answer for any maximal arity. The reason for thinking so has to do with
the nature of some proofs in this article and in~\cite{Kop17a}.
Moreover, there are related results where the arity 4 is essential, such as ~\cite[Corollaries~5.3--5.4]{KopBin}
and the theory of smoothly approximable structures~\cite[Theorems~2--6]{CH}.

The structure of this article is the following:
The next section recalls definitions and results which will be used later and clarifies the notation and terminology that will be used.
Section~\ref{Main results and their relationships}
states the main results and explains how they are related.
Once this has been done it only remains to prove 
Theorems~\ref{nontrivial dependence implies a nontrivial regular type},
\ref{nontrivial dependence implies nontrivial pregeometry}
and~\ref{impossibility of a structure satisfying the other result}.
The first two of these are proved in Section~\ref{Section about regular types}
and the third one in Section~\ref{impossibility of a nontrivial pregeometry}.

\section{Preliminaries}

\noindent
Familiarity with model theory including the basics of simplicity theory is assumed.
Elementary model theoretic results of relevance here can be found in \cite{Hod, TZ}.
For basic notions and results about simple theories we refer to any of the books \cite{Cas, Kim_book, Wag}.
The terminology and notation used here is relatively standard but nevertheless we make some clarifications.
A vocabulary (or signature) $V$ is {\em finite and relational} if it is finite and contains only relation symbols.
Its {\em maximal arity} is the maximal $k$ for which some symbol in it has arity $k$.
A finite relational vocabulary is called  {\em ternary} if its maximal arity is 3.
If $V$ is a ternary finite vocabulary then a $V$-structure is called {\em ternary}.
Structures will be denoted by calligraphic letters such as $\mcA, \ldots, \mcM, \mcN$ and their universes by
$A, \ldots, M, N$. Finite sequences/tuples of objects are denoted by $\bar{a}, \bar{b}, \ldots, \bar{x}, \bar{y}, \ldots$.
When writing $\bar{a} \in A$ it can mean that $\bar{a}$ belongs to $A$, but it can also mean that
each coordinate of the tuple $\bar{a}$ belongs to $A$; hopefully the meaning will be revealed by the context.
Occasionally we may write $\rng(\bar{a})$ for the set of all elements in the tuple $\bar{a}$ but we often abuse notation and
identify (when convenient) $\bar{a}$ with the {\em set} of elements in it.
Also, for sets $A$ and $B$ we often abbreviate `$A \cup B$' by `$AB$'.

Let $\mcM$ be a structure. Its complete theory is denoted by $Th(\mcM)$.
If $A \subseteq M$ then $\mcM \uhrc A$ denotes the substructure of $\mcM$ which is generated by $A$.
If $p$ is an $n$-type with respect to $\mcM$, then $p(\mcM)$ denotes the set of all $n$-tuples of elements from $M$ which realize $p$.
By $p \uhrc A$ we denote the restriction of $p$ to formulas with parameters from $A$.
Since the proofs will sometimes involve two structures with different complete theories (even different languages)
we will sometimes use the notation `$\ind^\mcM$', `$\tp_\mcM$', `$\acl_\mcM$' and `$\dcl_\mcM$' to clarify that
dividing/forking, types, algebraic closure and definable closure, respectively, is with respect to the structure $\mcM$.
By `$\tp^\mrqf_\mcM(\bar{a} / B)$' we denote the type of $\bar{a}$ over $B$ (in $\mcM$) 
{\em restricted to quantifier-free formulas}.
The notation `$\bar{a} \equiv_\mcM \bar{b}$' (`$\bar{a} \equiv^\mrqf_\mcM \bar{b}$')
means the same as `$\tp_\mcM(\bar{a}) = \tp_\mcM(\bar{b})$' (`$\tp^\mrqf_\mcM(\bar{a}) = \tp^\mrqf_\mcM(\bar{b})$').

For definitions of the notions {\em orthogonality} (of types), {\em regular types} and {\em weight} see \cite{Kim_book, Wag}.
A definition of a {\em pregeometry} is found in \cite[Definition~4.4.3]{Kim_book} and in several other books on model theory.
Recall also that in a simple theory {\em forking} and {\em dividing} are equivalent.

\begin{defin}\label{definition of minimal dependent set}{\rm
Let $T$ be a simple theory, $\mcM \models T$, $\bar{a}_1, \ldots, \bar{a}_n \in M\meq$, and $B \subseteq M\meq$.\\
(i) We call $A = \{\bar{a}_1, \ldots, \bar{a}_n\}$ {\em independent over $B$} if for every $\bar{a}_i \in A$,
$\bar{a}_i \underset{B}{\ind} (A \setminus \{\bar{a}_i\})$.\footnote{
A basic result is that $\{\bar{a}_1, \ldots, \bar{a}_n\}$ is independent over $B$ if and only if,
for every $0 < i < n$, $\bar{a}_{i+1} \underset{B}{\ind} \bar{a}_1, \ldots, \bar{a}_i$.}
Otherwise we call $A$ {\em dependent over $B$}. If $B = \es$ we may omit saying `over $B$'.\\
(ii) If $A = \{\bar{a}_1, \ldots, \bar{a}_n\}$ is dependent over $B$ and every proper subset of $A$ is
independent over $B$, then we call $A$ {\em minimal dependent over $B$} (and if $B = \es$ we may omit `over $B$').\\
(iii) A type $p$ over $B \subseteq M\meq$ is called {\em nontrivial} if $p(\mcM\meq)$ contains
a set of cardinality at least $3$ which is minimal dependent over $B$.
}\end{defin}

\noindent
Note that a minimal dependent set (over some set) must always be finite, by the finite character of dividing.
I have not been able to find the following useful result in the literature besides as an exercise in~\cite{TZ}, 
and its proof is indeed an elementary exercise in using dividing/forking.

\begin{fact}\label{useful fact about independence}{\rm (\cite[Exercise~7.2.6]{TZ})}
Let $T$ be simple, $\mcM \models T$, $\bar{a}_1, \ldots, \bar{a}_n \in M\meq$, $A \subseteq B \subseteq M\meq$
and suppose that $\bar{a}_1 \ldots \bar{a}_n \underset{A}{\ind} B$.
Then $\{\bar{a}_1, \ldots, \bar{a}_n\}$ is independent over $A$ if and only if it is independent over $B$.
\end{fact}

\begin{defin}\label{definition of minimal dependent set in a pregeometry}{\rm
Let $(X, \cl)$ be a pregeometry. \\
(i) We call $Y \subseteq X$ {\em independent} if, for every $y \in Y$, $y \notin \cl(Y \setminus \{y\})$.\footnote{
A basic result is that $Y = \{y_1, \ldots, y_n\}$ is independent if and only if, for every $0 < i < n$,
$y_{i+1} \notin \cl(y_1, \ldots, y_i)$.}
Otherwise we call $Y$ {\em dependent}.\\
(ii) If $Y \subseteq X$ is dependent and every proper subset of $Y$ is independent, then we call $Y$ {\em minimal dependent}.\\
(iii) Let $Y \subseteq X$. An independent set $Z \subseteq Y$ such that $Y \subseteq \cl(Z)$ is called a {\em basis} of $Y$ 
and $|Z|$ is called the {\em dimension} of $Y$; an elementary result is that all bases of $Y$ have the same cardinality 
so the dimension of $Y$ is well defined.\\
(iv) $(X, \cl)$ is {\em nontrivial} if there is $Y \subseteq X$ with cardinality at least $3$ such that $Y$ is minimal dependent.
}\end{defin}

\begin{fact}\label{existence of regular types} 
{\rm (\cite[Proposition~4.4.5]{Kim_book}, \cite[Proposition~5.1.11]{Wag})}
Suppose that $T$ is supersimple, $\mcM \models T$ and $A \subseteq M\meq$.
Let $n < \omega$ and let $p$ be a nonalgebraic $n$-type over $A$ which is realized in $M$.
Then, assuming that $\mcM$ is sufficiently saturated, there are $B \subseteq M\meq$ and a regular 1-type $q$
over $B$ such that $q$ is realized by a real element in $M$, $q$ is nonorthogonal to $p$ and $\su(q) \leq \su(p)$.
\end{fact}

\noindent
Some comments about the above fact may be in order.
Given $p$ as in the fact and a tuple $\bar{a} = (a_1, \ldots, a_n) \in M^n$ which realizes $p$,
$p$ will be nonorthogonal to $q' = \tp(a_1 /A)$ where (by the Lascar inequalities \cite{Kim_book, Wag}) $\su(q') \leq \su(q)$.
Now we can let $q$ be a 1-type of minimal SU-rank among all 1-types that are nonorthogonal to $p$ 
and realized by a real element, so $\su(q) \leq \su(p)$, and
then we can argue as in the proof of \cite[Proposition~4.4.5]{Kim_book} or \cite[Proposition~5.1.11]{Wag}
to show that $q$ is regular.

\begin{fact}\label{regular types and nondividing extensions}
{\rm (\cite[Remarks 4.4.2 (2) and~4.4.4]{Kim_book})}
Let $T$ be simple, $\mcM \models T$, $\bar{a} \in M\meq$ and $A \subseteq B \subseteq M\meq$. \\
(i) Suppose that $\bar{a} \underset{A}{\ind} B$. 
Then $\tp(\bar{a} / A)$ is regular if and only if $\tp(\bar{a} / B)$ is regular.\\
(ii) Suppose that $p$ is a regular type over $A$. Then $(p(\mcM), \cl)$ is a pregeometry if, for all
$b \in p(\mcM)$ and all $C \subseteq p(\mcM)$, $b \in \cl(C)$ if and only if $b \underset{A}{\nind} C$.
\end{fact}

\noindent
By the finite character of dividing the following conditions are equivalent when $T$ is simple:
\begin{itemize}
\item[(a)] For every $\mcM \models T$, all $A, B \subseteq M\meq$ and every finite tuple $\bar{a} \in M\meq$,
if $\bar{a} \underset{A}{\nind} B$ then there is $b \in B$ such that $\bar{a} \underset{A}{\nind} b$.
\item[(b)] For every $\mcM \models T$, all finite tuples $\bar{a}, \bar{b}, \bar{c} \in M\meq$ and every $A \subseteq M\meq$,
if $\bar{a} \underset{A}{\nind} \bar{b}\bar{c}$ then $\bar{a} \underset{A}{\nind} \bar{b}$ or 
$\bar{a} \underset{A}{\nind} \bar{c}$.
\end{itemize}

\begin{defin}\label{definition of trivial dependence}{\rm
A simple theory $T$ has {\em trivial dependence} if the two equivalent conditions~(a) and~(b) above hold.
(Otherwise it has {\em nontrivial dependence}.)
}\end{defin}

\noindent
Goode~\cite{Goode} studied a few variations of the notion of trivial dependence in the context of stable theories.
When he says `$T$ is totally trivial' it means the same as when we say `$T$ has trivial dependence'.
When he says `$T$ is trivial' it means that whenever $\bar{a}, \bar{b}$ and $\bar{c}$ are pairwise 
independent tuples over some set of elements, then $\{\bar{a}, \bar{b}, \bar{c}\}$ is independent
over the same set.
The next fact is stated in~\cite{Goode} for stable theories, but the proof uses only basic properties of forking/dividing which hold
also for simple structures.

\begin{fact}\label{trivial dependence for real elements gives trivial dependence}
{\rm (Goode \cite[Lemma 4]{Goode})}
Let $T$ be a simple theory and suppose that 
$\mcM \models T$, $A \subseteq M\meq$, $\bar{a}_1, \bar{a}_2, \bar{a}_3 \in M\meq$,
$\bar{a}_1 \underset{A}{\nind} \bar{a}_2 \bar{a}_3$, $\bar{a}_1 \underset{A}{\ind} \bar{a}_2$ and
$\bar{a}_1 \underset{A}{\ind} \bar{a}_3$.
Then there are tuples of real elements $\bar{b}_1, \bar{b}_2, \bar{b}_3 \in M$ such that 
$\bar{b}_1 \underset{A}{\nind} \bar{b}_2 \bar{b}_3$, $\bar{b}_1 \underset{A}{\ind} \bar{b}_2$ and
$\bar{b}_1 \underset{A}{\ind} \bar{b}_3$.
\end{fact}

\begin{defin}\label{definition of constraint}{\rm
Let $V$ be a vocabulary and $\mcM$ a $V$-structure.
A finite $V$-structure $\mcA$ is called a {\em constraint of $\mcM$} if $\mcA$ cannot be embedded into $\mcM$
but every proper substructure of $\mcA$ can be embedded into $\mcM$.
We say that $\mcM$ is {\em finitely constrained} if $\mcM$ has (up to isomorphism) only finitely many constraints.
}\end{defin}

\noindent
Every supersimple theory has elimination of hyperimaginaries 
(see \cite[Theorem~5.4.9]{Kim_book} or \cite[Theorem 5.3.1]{Wag}).
Also, every simple and `small' theory has elimination of hyperimaginaries
(see for example \cite[Corollary~5.3.5]{Kim_book}) and we note that $\omega$-categorical theories are small.
It follows that in the context of this article we need not consider hyperimaginary elements and 
two tuples of imaginaries have the same Lascar strong type (over some set) if and only if 
they have the same strong type (over the same set). For small, and hence for $\omega$-categorical,
simple theories this is perhaps most clearly stated
in \cite[Theorem~23]{Kim98}.
It follows that when we use the
independence theorem of simple theories (also called the type amalgamation theorem) \cite{Cas, Kim_book, Wag} 
we do not need to consider hyperimaginary elements; it suffices that the types that are to be amalgamated 
extend the same strong type over the ``base set'' rather than the same Lascar strong type.
When saying that a supersimple theory $T$ has {\em finite SU-rank} we mean that the SU-rank of every type with 
finitely many variables (realized by real elements) is finite.
The well known characterization of $\omega$-categorical theories by Engeler, Ryll-Nardzewski and Svenonius
(see \cite[Theorem~7.3.1]{Hod} or \cite[Theorem~4.3.1]{TZ} for example)
has the following consequence which will be used: if $T$ is $\omega$-categorical then every model of it is 
$\omega$-saturated and if $A$ is a {\em finite} subset of some model of $T$,
then there are only finitely many $n$-types over $A$ and each such type is isolated.

\section{Main results and their relationships}\label{Main results and their relationships}

\begin{theor}\label{nontrivial dependence implies a nontrivial regular type}
Suppose that $T$ is an $\omega$-categorical supersimple theory with nontrivial dependence. 
Then there are $\mcM \models T$, a finite set $C \subseteq M$ and a nontrivial regular 1-type $p$ over $C$ realized by real elements.
(Consequently, $T$ is not trivial in the sense of Goode~\cite{Goode} or Palac\'{i}n~\cite{Pal17}.)
\end{theor}

\noindent
If $p$ is a regular type like in Theorem~\ref{nontrivial dependence implies a nontrivial regular type},
then (by Fact~\ref{regular types and nondividing extensions}) $(p(\mcM), \cl)$ is a pregeometry if, 
for all $a \in p(\mcM)$ and all $A \subseteq p(\mcM)$,
$a \in \cl(A)$ if and only if $a \underset{C}{\nind} A$.
By the $\omega$-categoricity of $T$, for every $0 < n < \omega$, the relation $x_n \in \cl(x_1, \ldots, x_{n-1})$ on $p(\mcM)$
is $C$-definable.
Theorem~\ref{nontrivial dependence implies a nontrivial regular type} is proved
in Section~\ref{Section about regular types}. 
By continuing the proof of Theorem~\ref{nontrivial dependence implies a nontrivial regular type}, in the same section
with some extra assumptions added, we get the following:

\begin{theor}\label{nontrivial dependence implies nontrivial pregeometry}
Suppose that $V$ is a finite relational vocabulary and $T$ a complete $V$-theory such that
$T$ is supersimple with elimination of quantifiers. 
If $T$ has nontrivial dependence, then there is a finite relational vocabulary $V'$ with the same maximal arity as $V$ 
and a $V'$-structure $\mcM'$ satisfying the following conditions:
\begin{itemize}
\item[(i)] $Th(\mcM')$ has elimination of quantifiers and is supersimple.

\item[(ii)] All elements of $M'$ have the same type over $\acl_{(\mcM')\meq}(\es)$.

\item[(iii)] $(M', \cl)$ is a pregeometry if, for every $X \subseteq M'$ and every $x \in M'$, $x \in \cl(X)$ if and only if $x \nind^{\mcM'} X$.

\item[(iv)] $M'$ has a minimal dependent subset of cardinality at least 3
(where dependence is with respect to $Th(\mcM')$, or equivalently, with respect to the pregeometry from condition~(iii)).
\end{itemize}
\end{theor}

\noindent
In Section~\ref{impossibility of a nontrivial pregeometry} we prove the following:

\begin{theor}\label{impossibility of a structure satisfying the other result}
Let $V'$ be a ternary finite relational vocabulary. Then there does {\em not} exist a $V'$-structure $\mcM'$ such that (i)~--~(iv)
of Theorem~\ref{nontrivial dependence implies nontrivial pregeometry} hold.
(This holds also if we replace `supersimple' with `simple' in~(i).)
\end{theor}

\noindent
Recall the definition of {\em homogenous} structure from the introduction (which implies that a homogeneous structure
is countable and has a finite relational vocabulary).
By combining Theorems~\ref{nontrivial dependence implies nontrivial pregeometry}
and~\ref{impossibility of a structure satisfying the other result} and
Proposition~\ref{triviality of dependence implies finite rank}
below, we immediately get:

\begin{cor}\label{trivial dependence and finite rank}
If $V$ is a ternary finite relational vocabulary and $T$ a complete $V$-theory such that
$T$ is supersimple with elimination of quantifiers (so its unique countable model is homogeneous),
then $T$ has trivial dependence and finite SU-rank.
\end{cor}

\noindent
Corollary~\ref{trivial dependence and finite rank} easily gives the following improvement of Theorem~4.1 in~\cite{Kop17a}:

\begin{cor}\label{ternary homogeneous finitely constrained theories have trivial dependence}
Suppose that $\mcM$ is a ternary homogeneous finitely constrained simple structure.
Then $Th(\mcM)$ has trivial dependence and finite SU-rank.
\end{cor}

\noindent 
{\bf Proof.}
Suppose that $\mcM$ is ternary, homogeneous, finitely constrained and simple.
By \cite[Theorem~4.1]{Kop17a}, $Th(\mcM)$ is supersimple with finite SU-rank.
Then Corollary~\ref{trivial dependence and finite rank} implies that $Th(\mcM)$  has trivial dependence.
\hfill $\square$
\\

\noindent
The notion of trivial dependence used in this article implies the notion of `triviality' (or `1-triviality') used in Palac\'{i}n's article \cite{Pal17}
(which uses the terminology of Goode~\cite{Goode}).
Hence the next result was proved by Palac\'{i}n in that article.
Here we offer a different proof which is rather short, straightforward and avoids some technical notions such as
Lascar strong types, canonical bases, orthogonality and regular types.

\begin{prop}\label{triviality of dependence implies finite rank}{\rm (Palac\'{i}n \cite[Corollary 3.12]{Pal17})}
Suppose that $T$ is an $\omega$-categorical supersimple theory with trivial dependence.
Then $T$ has finite SU-rank.
\end{prop}

\noindent
{\bf Proof.}
Let $T$ be $\omega$-categorical and supersimple with trivial dependence.
Suppose that $T$ does not have finite SU-rank,
so there is $p \in S_1(T)$ with $\su(p) \geq \omega$.
Let $\mcM \models T$ be sufficently saturated so that all elements and sets that we talk about can be assumed
to come from $\mcM$.
Since $\su(p) \geq \omega$ there are (by for example~\cite[Remark~13.14]{Cas}) 
$A \subseteq M$ and an extension $p' \in S_1(A)$ of $p$ such that $\su(p') = \omega$.
Since $T$ is supersimple there is finite $A_0 \subseteq A$ such that $p'$ does not divide over $A_0$.
Let $p_0 = p' \uhrc A_0$, so $\su(p_0) = \su(p') = \omega$.
For every $0 < n < \omega$ there are $B_n \subseteq M$ and an
extension $p_n \in S_1(A_0B_n)$ of $p_0$ such that $\su(p_n) = n$.
Since $T$ has trivial dependence there is, for each $0 < n < \omega$, $b_n \in B_n$ such that 
$q_n = p_n \uhrc A_0b_n$ is a dividing extension of $p_0$. Then
\[
n = \su(p_n) \leq \su(q_n) < \su(p_0) = \omega.
\]
Let $k_1 = \su(q_1)$. 
Suppose that $0 < k_1 < k_2 < \ldots < k_m < \omega$ have been defined and that for each $i$, $k_i = \su(q_n)$ for some $n$.
Then choose $n < \omega$ such that $k_m < n$ and let $k_{m+1} = \su(q_n)$, so $k_{m+1} > k_m$.
Thus, by renaming (some of) the types and elements, we now have an infinite sequence
$0 < k_1 < k_2 < k_3 < \ldots$ of natural numbers, an infinite sequence of elemens $c_1, c_2, c_3, \ldots$ and, for each 
$0 < i < \omega$ a type $r_i \in S_1(A_0c_i)$ such that $\su(r_i) = k_i$.
Let $\bar{a}$ enumerate $A_0$ and, for every $0 < i < \omega$, let $d_i$ realize $r_i$.
Since
\[
\su(d_i / A_0c_i) \neq \su(d_j / A_0c_j) \ \text{ if } \ i \neq j
\]
it follows that $\tp(d_i, \bar{a}, c_i) \neq \tp(d_j, \bar{a}, c_j)$ if $i \neq j$.
Thus we have infinitely many $n$-types over $\es$ if $n = |\bar{a}| + 2$ and this contradicts that $T$ is $\omega$-categorical.
\hfill $\square$

\section{Proofs of Theorems~\ref{nontrivial dependence implies a nontrivial regular type}
and~\ref{nontrivial dependence implies nontrivial pregeometry}: 
nontrivial dependence gives (in the context) a nontrivial pregeometry}\label{Section about regular types}

\noindent
Let $T$ be an $\omega$-categorical supersimple theory with nontrivial dependence.
We also assume that $\mcM$ is a sufficiently saturated model of $T$, so that all elements or sets that we claim exist
(in some elementary extension of $\mcM\meq$) actually exist in $\mcM\meq$.
In this section, $\ind$, $\tp$, $\acl$ and $\dcl$ mean 
$\ind^{\mcM\meq}$, $\tp_{\mcM\meq}$, $\acl_{\mcM\meq}$ and $\dcl_{\mcM\meq}$, respectively.
Since $T$ has nontrivial dependence it follows from
Fact~\ref{trivial dependence for real elements gives trivial dependence}
that there are finite tuples $\bar{a}_1, \bar{a}_2, \bar{a}_3 \in M$ and $A \subseteq M\meq$ such that,
for some permutation $(i, j, k)$ of $(1, 2, 3)$, $\bar{a}_i \underset{A}{\nind} \bar{a}_j \bar{a}_k$,
$\bar{a}_i \underset{A}{\ind} \bar{a}_j$ and $\bar{a}_i \underset{A}{\ind} \bar{a}_k$.
Since $T$ is supersimple we may assume that 
\begin{itemize}
\item[(A)] there do not exist $\bar{a}'_i, \bar{a}'_j, \bar{a}'_k \in M$ and $A' \subseteq M\meq$ such that
\begin{itemize}
\item[] $\bar{a}'_i \underset{A'}{\nind} \bar{a}'_j \bar{a}'_k$,
$\bar{a}'_i \underset{A'}{\ind} \bar{a}'_j$, $\bar{a}'_i \underset{A'}{\ind} \bar{a}'_k$,
\item[] $\su(\bar{a}'_n / A') \leq \su(\bar{a}_n / A)$ for all $n = 1, 2, 3$, and
\item[] for some $1  \leq n \leq 3$, $\su(\bar{a}'_n / A') < \su(\bar{a}_n / A)$.
\end{itemize}
\end{itemize}

\noindent
Theorems~\ref{nontrivial dependence implies a nontrivial regular type}
and~\ref{nontrivial dependence implies nontrivial pregeometry}
will be proved via a sequence of lemmas. 
The proof of Theorem~\ref{nontrivial dependence implies a nontrivial regular type}
is finished after Lemma~\ref{finding a minimal dependent set in the same type}
and then the argument goes on, with added assumptions about $T$, to prove 
Theorem~\ref{nontrivial dependence implies nontrivial pregeometry}.

\begin{lem}\label{getting regular types}
There are real elements $a'_1, a'_2, a'_3 \in M$ and $B \subseteq M\meq$ such that
\begin{itemize}
\item[] $a'_i \underset{B}{\nind} a'_ja'_k$, $a'_i \underset{B}{\ind} a'_j$, $a'_i \underset{B}{\ind} a'_k$,
\item[] $\tp(a'_n /B)$ is regular for every $1 \leq n \leq 3$, and
\item[] $\su(a'_n / B) = \su(\bar{a}_n / A)$ for every $1 \leq n \leq 3$.
\end{itemize}
\end{lem}

\noindent
{\bf Proof.}
Let $p = \tp(\bar{a}_1 / A)$.
By 
Fact~\ref{existence of regular types} 
there is $A'' \subseteq M\meq$ and a regular type $p'' \in S_1^\mcM(A'')$ realized by a real element 
such that $p''$ and $p$ are nonorthogonal and $\su(p'') \leq \su(p)$.
This means that there are $B''$, $a''_1, \bar{a}^*_1$ such that
$A'', A \subseteq B''$,
$a''_1$ realizes a (complete)
nonforking extension of $p''$ to $B''$,
$\bar{a}^*_1$ realizes a nonforking extension of $p$ to $B''$ and
$a''_1 \underset{B''}{\nind} \bar{a}^*_1$.
Since $\tp(\bar{a}^*_1 /A) = \tp(\bar{a}_1 / A)$
there are $B \supseteq A$, $A' \subseteq B$ and $a'_1$ such that $p' = \tp(a'_1 / A')$ is regular,
$\su(p') = \su(p'')$ ($\leq \su(p)$) and
\begin{equation}\label{properties of a'-1 and bar-a-1}
\bar{a}_1 \underset{A}{\ind} B, \quad a'_1 \underset{A'}{\ind} B, \  \text{ and } \ \bar{a}_1 \underset{B}{\nind} a'_1.
\end{equation}
Since $a'_1 \underset{A'}{\ind} B$ it follows from 
Fact~\ref{regular types and nondividing extensions} 
that $\tp(a'_1 / B)$ is regular and $\su(a'_1 / B) = \su(p') \leq \su(p)$.
By the existence of nonforking extensions we may assume that 
\begin{equation}\label{a'-1B independent over}
a'_1B \underset{\bar{a}_1 A}{\ind} \bar{a}_2 \bar{a}_3.
\end{equation}

From~(\ref{properties of a'-1 and bar-a-1}) we get  $\bar{a}_1 \underset{A}{\ind} B$, so~(\ref{a'-1B independent over})
and transitivity gives
\begin{equation}\label{B independent from ... over A}
\bar{a}_1 \bar{a}_2 \bar{a}_3 \underset{A}{\ind} B.
\end{equation}
This together with 
Fact~\ref{useful fact about independence} 
implies that 
\begin{equation}\label{every proper subset is independent over B}
\text{for all $n, m \in \{1, 2, 3\}$, \ $\bar{a}_n \underset{B}{\ind} \bar{a}_m$ if and only if $\bar{a}_n \underset{A}{\ind} \bar{a}_m$.}
\end{equation}
Note that from~(\ref{B independent from ... over A}) it follows that $\su(\bar{a}_n / B) = \su(\bar{a}_n / A)$ for all $1 \leq n \leq 3$.
We now prove two claims and then explain how the lemma follows from these claims.

\bigskip

\noindent
{\bf Claim 1.}
Suppose that $\bar{a}_1 \underset{A}{\nind} \bar{a}_2 \bar{a}_3$, $\bar{a}_1 \underset{A}{\ind} \bar{a}_2$ and
$\bar{a}_1 \underset{A}{\ind} \bar{a}_3$.
Then $a'_1 \underset{B}{\nind} \bar{a}_2 \bar{a}_3$, $a'_1 \underset{B}{\ind} \bar{a}_2$ and $a'_1 \underset{B}{\ind} \bar{a}_2$.

\medskip

\noindent
{\bf Proof of Claim 1.}
By assumption, $\bar{a}_1 \underset{A}{\ind} \bar{a}_n$, for $n = 2, 3$, so~(\ref{a'-1B independent over})
and transitivity implies that $a'_1\bar{a}_1B \underset{A}{\ind} \bar{a}_n$ and hence
\begin{equation}\label{a'-1 ind from a-i}
a'_1 \underset{B}{\ind} \bar{a}_n \ \ \text{ and } \ \ \bar{a}_1 \underset{a'_1 B}{\ind} \bar{a}_n \ \ \text{ for } n = 2, 3.
\end{equation}
Towards a contradiction, suppose that $a'_1 \underset{B}{\ind} \bar{a}_2 \bar{a}_3$.
By assumption, $\bar{a}_1 \underset{A}{\nind} \bar{a}_2 \bar{a}_3$, so 
$a'_1 \bar{a}_1 B \underset{A}{\nind} \bar{a}_2 \bar{a}_3$. 
Transitivity and~(\ref{B independent from ... over A}) now implies that 
$a'_1 \bar{a}_1 \underset{B}{\nind} \bar{a}_2 \bar{a}_3$.
By assumption, $a'_1 \underset{B}{\ind} \bar{a}_2 \bar{a}_3$, so transitivity gives
$\bar{a}_1 \underset{a'_1 B}{\nind} \bar{a}_2 \bar{a}_3$.
From~(\ref{properties of a'-1 and bar-a-1}) it follows that
$\su(\bar{a}_1 / a'_1 B) < \su(\bar{a}_1 / B) = \su(\bar{a}_1 / A)$.
But now, taking $A' = a'_1B$, we have a situation which contradicts the assumption~(A).
Hence $a'_1 \underset{B}{\nind} \bar{a}_2 \bar{a}_3$, so the claim is proved.
\hfill $\square$

\bigskip

\noindent
{\bf Claim 2.}
Suppose that $\bar{a}_1 \bar{a}_2 \underset{A}{\nind} \bar{a}_3$, $\bar{a}_1 \underset{A}{\ind} \bar{a}_3$ and
$\bar{a}_2 \underset{A}{\ind} \bar{a}_3$.
Then $a'_1 \bar{a}_2 \underset{B}{\nind} \bar{a}_3$, $a'_1 \underset{B}{\ind} \bar{a}_3$ and $\bar{a}_2 \underset{B}{\ind} \bar{a}_3$.
The claim also holds if `$\bar{a}_2$' and `$\bar{a}_3$' switch places by letting `$\bar{a}_2$' and
`$\bar{a}_3$' switch roles in the proof.

\medskip

\noindent
{\bf Proof of Claim 2.}
By assumption, $\bar{a}_2 \underset{A}{\ind} \bar{a}_3$ so~(\ref{every proper subset is independent over B}) gives
$\bar{a}_2 \underset{B}{\ind} \bar{a}_3$.
The assumption that $\bar{a}_1 \underset{A}{\ind} \bar{a}_3$ together with~(\ref{a'-1B independent over})
and transitivity gives $a'_1 \bar{a}_1 B \underset{A}{\ind} \bar{a}_3$ and hence $a'_1 \underset{B}{\ind} \bar{a}_3$.

Towards a contradiction, suppose that $a'_1 \bar{a}_2 \underset{B}{\ind} \bar{a}_3$.
By assumption, $\bar{a}_1 \bar{a}_2 \underset{A}{\nind} \bar{a}_3$ so 
\[
a'_1 \bar{a}_1 \bar{a}_2 B \underset{A}{\nind} \bar{a}_3.
\]
By assumption, $\bar{a}_1 \underset{A}{\ind} \bar{a}_3$ and by~(\ref{a'-1B independent over})
we have $a'_1 \bar{a}_1 B \underset{A}{\ind} \bar{a}_3$, so $\bar{a}_1 \underset{a'_1 B}{\ind} \bar{a}_3$.
The assumption that $a'_1 \bar{a}_2 \underset{B}{\ind} \bar{a}_3$ gives $\bar{a}_2 \underset{a'_1 B}{\ind} \bar{a}_3$.
From~(\ref{properties of a'-1 and bar-a-1})
we get $\su(\bar{a}_1 / a'_1 B) < \su(\bar{a}_1 / B) = \su(\bar{a}_1 / A)$.
Thus we have, with $A' = a'_1 B$, a situation that contradicts the assumption~(A).
Hence $a'_1 \bar{a}_2 \underset{B}{\nind} \bar{a}_3$, so the claim is proved.
\hfill $\square$
\\

\noindent
Let $a^*_1 = a'_1$, $a^*_2 = \bar{a}_2$ and $a^*_3 = \bar{a}_3$.
If $i = 1$ then, by Claim~1, $a^*_i \underset{B}{\nind} a^*_j a^*_k$, 
$a^*_i \underset{B}{\ind} a^*_j$ and $a^*_i \underset{B}{\ind} a^*_k$.
If $i = 2$ or $i = 3$ then we use Claim~2 to get the same conclusion.
Next, starting with $a'_1, \bar{a}_2, \bar{a}_3$ and $B$ instead of 
$\bar{a}_1, \bar{a}_2, \bar{a}_3$ and $A$, 
we argue in the same way as before Claim~1 to get $B' \supseteq B$ and $a'_2$ such that 
$\tp(a'_2 / B')$ is regular, $\su(a'_2 / B') \leq \su(\bar{a}_2 / B)$,
(\ref{properties of a'-1 and bar-a-1}) holds if $\bar{a}_1$, $a'_1$, $A$ and $B$ are replaced by $\bar{a}_2$, $a'_2$, $B$ and $B'$, respectively,
$a'_2 B' \underset{\bar{a}_2 B}{\ind} a'_1\bar{a}_3$ (which corresponds to~(\ref{a'-1B independent over}))
and~(\ref{B independent from ... over A})--(\ref{every proper subset is independent over B}) hold if 
$\bar{a}_1$, $A$ and $B$ are replaced by $a'_1$, $B$ and $B'$, respectively.
Then Claims~1 and~2 hold if $\bar{a}_1, \bar{a}_2, \bar{a}_3$, $a'_1$, $A$ and $B$ are replaced by
$\bar{a}_2, a'_1, \bar{a}_3$, $a'_2$, $B$ and $B'$, respectively (so $\bar{a}_2$ and $a'_2$ take the ``active'' role now).
Now let $a^*_1 = a'_1$, $a^*_2 = a'_2$ and $a^*_3 = \bar{a}_3$.
If $i = 2$ then, by Claim~1, $a^*_i \underset{B'}{\nind} a^*_j a^*_k$, 
$a^*_i \underset{B'}{\ind} a^*_j$ and $a^*_i \underset{B'}{\ind} a^*_k$.
Otherwise we get the same conclusion by Claim~2.
Now we repeat the same kind of argument a final round to get $a'_3$ and $B'' \supseteq B$ such that
$\tp(a'_3 / B'')$ is regular, $\su(a'_3 / B'') \leq \su(\bar{a}_3 / B')$,
$a'_i \underset{B''}{\nind} a'_j a'_k$, 
$a'_i \underset{B''}{\ind} a'_j$ and $a'_i \underset{B''}{\ind} a'_k$.
For every $1 \leq n \leq 3$ we have $\su(a'_n / B'') \leq \su(\bar{a}_n / A)$
(this follows from the three versions of~(\ref{B independent from ... over A}) that have been obtained during the three ``rounds'').
Hence the assumption~(A) implies that $\su(a'_n / B'') = \su(\bar{a}_n / A)$ for every $1 \leq n \leq 3$.
So if $B''$ is renamed as $B$, then we have proved the statement of the lemma.
\hfill $\square$
\\

\noindent
Let $a'_1, a'_2, a'_3$ and $B$ be as in Lemma~\ref{getting regular types} and to
simplify notation we assume that $(i, j, k) = (1, 2, 3)$.

\begin{lem}\label{a finite base set of reals}
There is a finite set (of real elements) $C \subseteq M$ such that 
\begin{itemize}
\item[(i)] $\tp(a'_n / C)$ is regular for every $1 \leq n \leq 3$,
\item[(ii)] $\su(a'_n / C) = \su(a'_n / B)$ for every $1 \leq n \leq 3$, and 
\item[(iii)] $a'_1 \underset{C}{\nind} a'_2 a'_3$ and $a'_1 \underset{C}{\ind} a'_n$ for every $2 \leq n \leq 3$.
\end{itemize}
\end{lem}

\noindent
{\bf Proof.}
Let $B' \subseteq M$ be such that $B \subseteq \dcl(B')$.
We may assume that $B' \underset{B}{\ind} a'_1 a'_2 a'_3$, so in particular $a'_n \underset{B}{\ind} B'$ for each $n$.
Since $B \subseteq \acl(B')$ it follows from
Fact~\ref{useful fact about independence}  
that $a'_1 \underset{B'}{\nind} a'_2 a'_3$ and $a'_1 \underset{B'}{\ind} a'_n$ for $n = 2, 3$.
By Fact~\ref{regular types and nondividing extensions}, 
$\tp(a'_n / B')$ is regular for each $n = 1, 2, 3$.

Since $T$ is supersimple there is finite $C \subseteq B'$ such that
$a'_1 a'_2 a'_3 \underset{C}{\ind} B'$.
Again it follows from 
Facts~\ref{useful fact about independence} 
and~\ref{regular types and nondividing extensions}
that $a'_1 \underset{C}{\nind} a'_2 a'_3$, $a'_1 \underset{C}{\ind} a'_n$ for $n = 2, 3$
and $\tp(a'_n / C)$ is regular for each $n = 1, 2, 3$.
Moreover, we have that $\su(a'_n / C) = \su(a'_n / B') = \su(a'_n / B)$.
\hfill $\square$
\\

\noindent
Let $C$ be as in Lemma~\ref{a finite base set of reals} and, for $n = 1, 2, 3$, let $p_n = \tp(a'_n / C)$.

\begin{lem}\label{getting a minimal dependent set}
We have $a'_2 \underset{C}{\ind} a'_3$, so $\{a'_1, a'_2, a'_3\}$ is minimal dependent over $C$.
\end{lem}

\noindent
{\bf Proof.}
For every $1 \leq n \leq 3$, $p_n$ is regular so it has weight~1 
(see \cite[Proposition~4.4.9]{Kim_book} or \cite[Lemma~5.2.11]{Wag})
from which it follows that the relation 
$x \underset{C}{\nind} y$ is an equivalence relation on 
$X = p_1(\mcM) \cup p_2(\mcM) \cup p_3(\mcM)$. Denote the relation $x \underset{C}{\nind} y$ on $X$ by $E$.
Since $T$ is $\omega$-categorical and $C$ is a finite set of reals, $E$ is $C$-definable.
Let $\bar{c}$ enumerate $C$.
For $\bar{x}, \bar{y}$ of the same length as $\bar{c}$, define 
\[
\text{$F(\bar{x}x, \bar{y}y)$ if and only if
$\bar{x} = \bar{y}$, $\tp(\bar{x}) = \tp(\bar{c})$ and $x \underset{\bar{x}}{\nind} y$.}
\]
Then $F$ is an equivalence relation which is $\es$-definable, so every $F$-class is represented by an imaginary element, and hence
every $E$-class is also represented by an imaginary element.
Let $b$ be the imaginary element representing the $E$-class of $a'_2$ (or strictly speaking the $F$-class of $\bar{c}a'_2$),
so $b \in \dcl(a'_2C)$.

Suppose  for a contradiction that $a'_2 \underset{C}{\nind} a'_3$.
If there would only be finitely many $E$-classes, then $b \in \acl(C)$ and hence $a'_2 \underset{Cb}{\nind} a'_3$.
By the existence of nonforking extensions we could find $a$ in the same $E$-class as $a'_2$ such that $a'_2 \underset{Cb}{\ind} a$,
which (as $b \in \acl(C)$) implies that $a'_2 \underset{C}{\ind} a$,
but this contradicts that (by definition) $E(x, y)$ implies $x \underset{C}{\nind} y$.
Hence there are infinitely many $E$-classes and thus $b \notin \acl(C)$.
From the assumption that $a'_2 \underset{C}{\nind} a'_3$ we get $E(a'_2, a'_3)$ and hence
$b \in \dcl(a'_3C)$. Consequently $a'_n \underset{C}{\nind} b$ for $n = 2, 3$ and hence
$\su(a'_n / bC) < \su(a'_n / C)$ for $n = 2, 3$.
By the choice of $a'_1, a'_2$ and $a'_3$, we have 
$a'_1 \underset{C}{\nind} a'_2 a'_3$ and $a'_1 \underset{C}{\ind} a'_n$ for $n = 2, 3$.
As $b \in \acl(a'_n C)$ for $n = 2, 3$
we get $a'_1 \underset{C}{\ind} a'_nb$ and hence
$a'_1 \underset{C}{\ind} b$ and $a'_1 \underset{bC}{\ind} a'_n$ for $n = 2, 3$.
Using that $a'_1 \underset{C}{\nind} a'_2 a'_3 b$ and transitivity it also follows that $a'_1 \underset{bC}{\nind} a'_2 a'_3$.
But since $\su(a'_n / C) = \su(\bar{a}_n / A)$ for all $n = 1, 2, 3$, 
this situation contradicts the assumption~(A).
\hfill $\square$

\begin{lem}\label{finding a minimal dependent set in the same type}
For some $1 \leq n \leq 3$, $p_n(\mcM)$ contains a set of cardinality at least 3 which is minimal dependent over $C$.
\end{lem}

\noindent
{\bf Proof.}
If all of $a'_1, a'_2$ and $a'_3$ realize the same $p_n$ then we are done.
So suppose that this is not the case.
Without loss of generality we can assume that $a'_2, a'_3 \notin p_1(\mcM)$.
By the existence of nonforking extensions there are $a''_2, a''_3 \in p_2(\mcM) \cup p_3(\mcM)$ so that 
$\tp(a''_2, a''_3 / a'_1 C) = \tp(a'_2, a'_3 / a'_1 C)$ and $a'_2 a'_3 \underset{a'_1 C}{\ind} a''_2 a''_3$.
By Lemma~\ref{getting a minimal dependent set},
$\{a'_1, a'_2, a'_3\}$ is minimal dependent over $C$ and hence
$a'_1 \underset{C}{\nind} a'_2 a'_3$ and consequently
$a'_1 \underset{C}{\nind} a''_2 a''_3$.
Since $p_1$ is regular and therefore has weight 1 it follows that 
$a'_2 a'_3 \underset{C}{\nind} a''_2 a''_3$ and consequently the set $\{a'_2, a'_3, a''_2, a''_3\}$ is dependent over $C$.
The fact that  $\{a'_1, a'_2, a'_3\}$ is minimal dependent over $C$
together with $a'_2 a'_3 \underset{a'_1 C}{\ind} a''_2 a''_3$ and transitivity implies that
every proper subset of $Y = \{a'_2, a'_3, a''_2, a''_3\}$ is independent over $C$.
Note also that $p_1$ is not realized by any element in $Y$.

If only one of $p_2$ and $p_3$ is realized in $Y$ then we are done.
So suppose that this is not the case.
Then, by the choice of the elements, two elements in $Y$ realize $p_2$ and two elements realize $p_3$.
So by renaming the elements we have $Y = \{b_1, b_2, b_3, b_4\}$ where $b_3$ and $b_4$ realize $p_2$.
Choose $b'_1, b'_2 \in p_3(\mcM)$ so that 
\begin{equation}\label{the copy of b-1 b-2 over b-3 and b-4}
\text{$\tp(b'_1, b'_2, b_3 / b_4 C) = \tp(b_1, b_2, b_3 / b_4 C)$  \ and \ 
$b_1 b_2 \underset{C b_3 b_4}{\ind} b'_1 b'_2$.}
\end{equation}
Let $n \in \{1, 2\}$. Since $\{b_n, b_3, b_4\}$ is independent over $C$ it follows from~(\ref{the copy of b-1 b-2 over b-3 and b-4})
and transitivity that
$b_n \underset{C}{\ind} b'_1 b'_2 b_3 b_4$.
From~(\ref{the copy of b-1 b-2 over b-3 and b-4}) it follows that  $\{b'_1, b'_2, b_3\}$ is independent over $C$
and hence $\{b_n, b'_1, b'_2, b_3\}$ is independent over $C$.
In the same way it follows that $\{b_1, b_2, b'_n, b_3\}$ is independent over $C$.

Since $\{b_3, b_4\}$ is independent over $C$ it follows that $b_4$ realizes a nonforking extension of $p_2$ to $Cb_3$.
As $Y$ is minimal dependent over $C$ we have $b_1 b_2 \underset{C b_3}{\nind} b_4$ and 
(using~(\ref{the copy of b-1 b-2 over b-3 and b-4}))
$b'_1 b'_2 \underset{C b_3}{\nind} b_4$.
Since $p_2$ has weight 1 it follows that $b_1 b_2 \underset{C b_3}{\nind} b'_1 b'_2$ and hence
$\{b_1, b_2, b'_1, b'_2, b_3\}$ is dependent over $C$.
If already $\{b_1, b_2, b'_1, b'_2\}$ is dependent over $C$, then it is minimal dependent and only $p_3$ is realized in it, so we are done.
If it is independent over $C$, then $Y' = \{b_1, b_2, b'_1, b'_2, b_3\}$ is minimal dependent over $C$ and all elements in $Y'$ except $b_3$
realizes $p_3$.
This means that we can argue similarly as we did in the beginning of the proof of this lemma to get 
a set of cardinality $8$ which is minimal dependent and in which all elements realize $p_3$.
\hfill $\square$
\\

\noindent
By Lemma~\ref{finding a minimal dependent set in the same type}, for some $1 \leq n \leq 3$,
$p_n(\mcM)$ contains a set of cardinality at least 3 which is minimal dependent over $C$.
This proves Theorem~\ref{nontrivial dependence implies a nontrivial regular type}.

\begin{rem}\label{remark about the use of omega-categoricity}{\rm 
In proving Theorem~\ref{nontrivial dependence implies a nontrivial regular type}
the assumption that $T$ is $\omega$-categorical was used only once
and it was in the proof of Lemma~\ref{getting a minimal dependent set}.
There $\omega$-categoricity was used to conclude that the relation $x \underset{C}{\nind} y$
on $p_1(\mcM) \cup p_2(\mcM) \cup p_3(\mcM)$ is $C$-definable.
In fact, in every simple theory, if $p \in S(C)$ then the set
$\{(a, b) : \text{ $a$ realizes $p$ and }  a \underset{C}{\ind} b\}$ is type-definable over $C$.
As a finite union of type-definable sets is type-definable, the relation $x \underset{C}{\ind} y$
is type-definable on $p_1(\mcM) \cup p_2(\mcM) \cup p_3(\mcM)$.
So if the relation $x \underset{C}{\nind} y$ is type-definable on $p_1(\mcM) \cup p_2(\mcM) \cup p_3(\mcM)$,
then, by compactness, it is definable.
However, if $T$ is not $\omega$-categorical then the relation $x \underset{C}{\nind} y$ need not be type-definable.
The example called $T_2$ which follows after the proof of Proposition~5 in~\cite{Goode}
is superstable (but not $\omega$-categorical)
and every regular type is trivial, but $T_2$ is not `totally trivial' in the sense of \cite{Goode} which means that
it has nontrivial dependence in the sense of this article.
Moreover, there is a regular 1-type $q$ of $T_2$ over $\es$ (with SU-rank $\omega$) such that the relation 
$x \nind y$ on $q$ is not type-definable.
}\end{rem}

\noindent
By Lemma~\ref{finding a minimal dependent set in the same type} we may, without loss of generality,
assume that $p_1(\mcM)$ contains a set of cardinality at least $3$ which is minimal dependent over $C$.
By 
Fact~\ref{regular types and nondividing extensions}, 
$(p_1(\mcM), \cl)$ is a pregeometry if $\cl$ is defined as follows for $x \in p(\mcM)$ and $X  \subseteq p_1(\mcM)$:
$x \in \cl(X)$ if and only if $x \underset{C}{\nind} X$.

We now continue with the proof of Theorem~\ref{nontrivial dependence implies nontrivial pregeometry}.
Thus, {\bf \em we now add the assumption that $T$ is a $V$-theory with elimination of quantifiers where $V$ is a finite
relational vocabulary.}
Let $P = \{p_1\}$.
(The choice of notation may seem awkward in the present context, but makes the correspondence to results in
\cite[Section~3]{Kop17a} clear.)
We now define a new vocabulary $V_C$, which will be finite and relational with the same maximal arity as $V$, 
and then we define a $V_C$-structure called $\mcM_P$.
(This corresponds to Definition~3.1 in~\cite{Kop17a}.)

\begin{defin}\label{definition of structure M-p generated by p}{\rm 
(i) Let $V_C$ be a finite relational vocabulary such that $V \subseteq V_C$ and for every $R \in V$ of arity $r > 1$, every 
$0 < k < r$, every permutation $\pi$ of $\{1, \ldots, r\}$, 
and every $\bar{a} \in C^k$, $V_C$ has a relation symbol $Q_{R, \bar{a}, \pi}$ of arity $r-k$. 
We also assume that $V_C$ has no other symbols than those described. 

\medskip
\noindent
(ii) Let $\mcM_P$ be the (infinite) $V_C$-structure with universe $M_P = p_1(\mcM)$ 
and where the symbols in $V_C$ are interpreted as follows:
\begin{itemize}
\item[(a)] If $R \in V$ has arity $r$, then $R^{\mcM_P} = R^\mcM \cap (M_P)^r$.

\item[(b)] If $Q_{R, \bar{a}, \pi} \in V_C \setminus V$ where $R \in V$ has arity $r$ and $|\bar{a}| = k$, 
then for every $\bar{b} \in (M_P)^{r-k}$,
$\bar{b} \in (Q_{R, \bar{a}, \pi})^{\mcM_P}$ if and only if $\pi(\bar{b}\bar{a}) \in R^\mcM$
(where, if $\bar{b}\bar{a} = (d_1, \ldots, d_r)$ then $\pi(\bar{b}\bar{a}) = (d_{\pi(1)}, \ldots, d_{\pi(r)})$).
\end{itemize}
}\end{defin}

\noindent
In~\cite[Lemmas 3.2--3.4]{Kop17a} the following was proved:

\begin{lem}\label{properties of M-P}
(i) Let $\bar{c}$ be an enumeration of $C$. For all $\bar{a}, \bar{b} \in M_P$,
\[
\bar{a} \equiv_{\mcM_P}^\mrqf \bar{b} \quad \text{ if and only if } \quad \bar{a}\bar{c} \equiv_\mcM^\mrqf \bar{b}\bar{c}.
\]
(ii) $Th(\mcM_P)$ has elimination of quantifiers and is simple.\footnote{
Lemma~3.2 in~\cite{Kop17a} uses the terminology `homogeneous' instead of `elimination of quantifiers'.
But $\mcM$ has elimination of quantifiers if and only if its unique, up to isomorphism, countable elementary substructure 
has elimination of quantifiers, and the same holds for $\mcM_P$. Therefore we can formulate part~(ii) as done here.}\\
(iii) For all $\bar{a}, \bar{b} \in M_P$, $\bar{a} \ind^{\mcM_P} \bar{b}$ if and only in $\bar{a} \underset{C}{\ind}^\mcM \bar{b}$.
\end{lem}

\begin{rem}\label{interpretability and supersimplicity}{\rm
(i) One characterization of supersimplicity is that for all $\bar{a}$ and $B$ there is a finite 
$B' \subseteq B$ such that the type of $\bar{a}$ over $B$ does not fork over $B'$.
Using this characterization it is easy to prove that if $\mcM$ is supersimple and $d_1, \ldots, d_n \in M$, 
then the expansion of $\mcM$ with constants for the elements $d_1, \ldots, d_n$ is also supersimple.\\
(ii) By \cite[Remark~2.8.14]{Wag}, every theory which is interpretable in a supersimple theory is supersimple.
It follows that $\mcM_P$ is supersimple.
From these observations it also follows that if $d_1, \ldots, d_n \in M\meq$, then
the expansion of $\mcM\meq$ by constants for the elements $d_1, \ldots, d_n$ is also supersimple and every 
structure which is interpretable in this expansion is supersimple.
}\end{rem}

\begin{defin}\label{definition of the equivalence relation}{\rm
(i) For $x \in M_P$ and $X \subseteq M_P$, define $x \in \cl(X)$ if and only if $x \underset{C}{\nind^\mcM} X$.\\
(ii) For $x, y \in M_P$, define $x \sim y$ if and only if 
$\tp(x / \acl(C)) = \tp(y / \acl(C))$, where types and algebraic closure are taken in $\mcM\meq$.\\
(iii) If $X$ is an equivalence class of `$\sim$', then define for all $x \in X$ and all $Y \subseteq X$,
$\cl_X(Y) = \cl(Y) \cap X$.
}\end{defin}

\noindent 
Note that `$\sim$' has only finitely many equivalence classes (as $C$ is finite and $T$ is $\omega$-categorical).
Moreover, since $p_1$ is nonalgebraic and regular every $\sim$-class is infinite.
By 
Fact~\ref{regular types and nondividing extensions}, 
$(M_P, \cl)$ is a pregeometry. 
Since $\mcM$ is $\omega$-categorical, $M_P$ is a $C$-definable set in $\mcM$.
Note also that, by part~(i) of Lemma~\ref{properties of M-P}, 
for every relation $R \subseteq (M_P)^n$, $R$ is $\es$-definable in $\mcM_P$ 
if and only if it is $C$-definable in $\mcM$.
By $\omega$-categoricity, $\sim$ is $C$-definable in $\mcM$ and hence 
(by Lemma~\ref{properties of M-P}~(i)) it is
$\es$-definable in $\mcM_P$.

\begin{lem}\label{the structure on a single equivalence class}
Let $X$ be any equivalence class of $\sim$ and let $\mcN  = \mcM_P \uhrc X$. Then:\\
(i) $Th(\mcN)$ has elimination of quantifiers and is supersimple.\\
(ii) For all $\bar{a}, \bar{b} \in X$, $\bar{a} \ind^{\mcN} \bar{b}$ if and only if $\bar{a} \ind^{\mcM_P} \bar{b}$.\\
(iii) All elements of $X$ have the same type, in $\mcN\meq$, over $\acl_{\mcN\meq}(\es)$. \\
(iv) $(X, \cl_X)$ is a pregeometry and for every $x \in X$ and every $Y \subseteq X$, $x \in \cl_X(Y)$ if and only if $x \nind^\mcN Y$.
\end{lem}

\noindent
{\bf Proof.}
(i) To show that $Th(\mcN)$ has elimination of quantifiers it suffices (by the use of back and forth arguments
or Ehrenfeucht-Fra\"{i}ss\'{e} games, see for example~\cite[Section~3.3]{Hod}) to show that if
$\bar{a}, \bar{b} \in X$, $\bar{a} \equiv_{\mcM_P}^\mrqf \bar{b}$ and $d \in X$, then there is $e \in X$
such that $\bar{a}d \equiv_{\mcM_P}^\mrqf \bar{b}e$. 
For this it suffices to show that if 
$\bar{a}, \bar{b} \in X$, $\bar{a}\bar{c} \equiv_\mcM^\mrqf \bar{b}\bar{c}$ and $d \in X$, then there is $e \in X$
such that $\bar{a}\bar{c}d \equiv_\mcM^\mrqf \bar{b}\bar{c}e$. 
But the later implication follows directly since $\mcM$ has elimination of quantifiers and $\sim$ is definable in $\mcM$
with parameters from $\bar{c}$.

That $Th(\mcN)$ is supersimple follows from Remark~\ref{interpretability and supersimplicity}
since $\mcN$ is interpretable in $\mcM\meq$ using the parameters in $C$ and the imaginary element
which corresponds to the equivalence class $X$.

(ii) Let $\bar{a}, \bar{b} \in X$ and let $\varphi(\bar{x}, \bar{y})$ be a quantifier free $V_C$-formula that 
isolates $\tp_\mcN(\bar{a}, \bar{b})$.
Also let $E(x, y)$ be a quantifier free $V_C$-formula which, in $\mcM_P$, defines `$\sim$'.
Then, whenever $x_i, x_j \in \rng(\bar{x})$ and $y_k, y_l \in \rng(\bar{y})$, we have
\begin{equation}\label{x-i and y-j in the same sim-class}
\models \forall \bar{x}, \bar{y} \Big(\varphi(\bar{x}, \bar{y}) \ \rightarrow \ 
\big(E(x_i, x_j) \wedge E(y_k, y_l) \wedge E(x_i, y_k) \big) \Big).
\end{equation}

First suppose that $\bar{a} \nind^\mcN \bar{b}$, so there are $\bar{b}_i \in X$, $i < \omega$, such that
$(\bar{b}_i : i < \omega)$ is an $\es$-indiscernible sequence in $\mcN$, $\bar{b}_0 = \bar{b}$ and
$\{\varphi(\bar{x}, \bar{b}_i) : i < \omega\}$ is $k$-inconsistent (with respect to $\mcN$) for some $k < \omega$.
As $\mcN$ is a substrructure of $\mcM_P$ and both have elimination of quantifiers, $(\bar{b}_i : i < \omega)$
is $\es$-indiscernible in $\mcM_P$ as well.
If there would be $\bar{a}' \in M_P$ such that $\mcM \models \bigwedge_{i=0}^{k-1} \varphi(\bar{a}', \bar{b}_i)$,
then, since some $a' \in \rng(\bar{a}')$ must not belong to $X$, we get, for any $b \in \rng(\bar{b}_0)$,
\[
\mcM_P \models \varphi(\bar{a}', \bar{b}_0) \wedge \neg E(a', b)
\]
and this contradicts~(\ref{x-i and y-j in the same sim-class}).
Hence $\{\varphi(\bar{x}, \bar{b}_i) : i < \omega\}$ is $k$-inconsistent with respect to $\mcM_P$
so $\bar{a} \nind^{\mcM_P} \bar{b}$.

Now suppose that $\bar{a} \nind^{\mcM_P} \bar{b}$.
Since $\mcN$ is a substructure of $\mcM_P$ and both have elimination of quantifiers it follows that
$\varphi(\bar{x}, \bar{y})$ also isolates $\tp_{\mcM_P}(\bar{a}, \bar{b})$.
Hence there are $\bar{b}_i \in M_P$ such that
$(\bar{b}_i : i < \omega)$ is $\es$-indiscernible,
$\bar{b}_0 = \bar{b}$ and $\{\varphi(\bar{x}, \bar{b}_i) : i < \omega\}$ is $k$-inconsistent 
(with respect to $\mcM_P$) for some $k < \omega$.
From~(\ref{x-i and y-j in the same sim-class}) it follows that, for all $i < \omega$, all elements in 
$\rng(\bar{b}_i)$ belong to the same $\sim$-class.
Since $(\bar{b}_i : i < \omega)$ is $\es$-indiscernible we have either
\begin{itemize}
\item for all $i < j$, all elements in $\rng(\bar{b}_i) \cup \rng(\bar{b}_j)$ belong to the same $\sim$-class, or
\item for all $i < j$, every $b \in \rng(\bar{b}_i)$ belongs to a different $\sim$-class than any $b' \in \rng(\bar{b}_j)$.
\end{itemize}
However, as there are only finitely many $\sim$-classes it follows that we are in the first case.
Since $\bar{b}_0 = \bar{b} \in X$ it follows that $\bar{b}_i \in X$ for all $i < \omega$.
If $\{\varphi(\bar{x}, \bar{b}_i) : i < \omega\}$ would be $k$-consistent with respect to $\mcN$, then, 
as $\mcN \subseteq \mcM_P$ and $\varphi$ is quantifier free, the same set would be $k$-consistent with 
respect to $\mcM_P$, which contradicts the assumption. Hence we conclude that $\bar{a} \nind^\mcN \bar{b}$.

(iii) Suppose, for a contradiction, that there are elements $a, b \in X$ such that 
\[
\tp_{\mcN\meq}(a / \acl_{\mcN\meq}(\es)) \neq
\tp_{\mcN\meq}(b / \acl_{\mcN\meq}(\es)).
\]
Then there is a nontrivial equivalence relation `$\approx$' on $X$
which is $\es$-definable in $\mcN$ and such that $a$ and $b$ belong to different $\approx$-classes. 
Define an equivalence relation on $M$ as follows:
\[
F(x, y) \ \Longleftrightarrow \ p_1(x) \ \wedge \ p_1(y) \ \wedge \ x \sim y \ \wedge \  x \approx y.
\]
Then $F$ is $C$-definable in $\mcM$.
Since $a, b \in X$ and since they belong to different $F$-classes it follows that
$a$ and $b$ have different types, in $\mcM\meq$, over $\acl_{\mcM\meq}(C)$ so they belong
to different $\sim$-classes which contradicts
that $a, b \in X$.

(iv) We already noted (after Definition~\ref{definition of the equivalence relation})
that $(M_P, \cl)$ is a pregeometry and from this it is a straightforward exercise to show that $(X, \cl_X)$ is a pregeometry.
Let $Y \subseteq X$ and $x \in X$.
Suppose that $x \in \cl_X(Y)$.
Then $x \in \cl(Y) \cap X$, so $x \underset{C}{\nind}^\mcM Y$ and hence $x \underset{C}{\nind}^\mcM \bar{y}$
for some finite tuple $\bar{y} \in Y$.
By Lemma~\ref{properties of M-P}~(iii), $x \nind^{\mcM_P} \bar{y}$ and by part~(ii) of this lemma,
$x \nind^\mcN \bar{y}$ so $x \nind^\mcN Y$.

Now suppose that $x \nind^\mcN Y$, so $x \nind^\mcN \bar{y}$ for some finite tuple $\bar{y} \in Y$.
By part~(ii) of this lemma and part~(iii) of Lemma~\ref{properties of M-P},
$x \underset{C}{\nind}^\mcM \bar{y}$, so $x \underset{C}{\nind}^\mcM Y$
and hence $x \in \cl(Y)$ (and by assumption $x \in X$).
\hfill $\square$

\begin{lem}\label{all elements have the same type}
There is an 
equivalence class $X$ of `$\sim$' such that $X$ contains a set of cardinality at least 3 which is minimal dependent over $\es$
with respect to $\mcM_P$.
\end{lem}

\noindent
{\bf Proof.}
We know that $M_P$ ($= p_1(\mcM)$) contains a set of cardinality at least 3 which is minimal dependent over $C$ when dividing
is considered with respect to $\mcM$.
By Lemma~\ref{properties of M-P}~(iii), this set is minimal dependent over $\es$ when dividing is considered with respect to $\mcM_P$.
Now Lemma~\ref{all elements have the same type} follows by using 
Lemma~\ref{reducing the number of types} below as many times as needed.
\hfill $\square$

\begin{lem}\label{reducing the number of types}
Suppose that $\mcN$ is a simple $\omega$-saturated structure and $(N, \cl_\mcN)$ a pregeometry such that
for every $x \in N$ and every $X \subseteq N$, $x \in \cl_\mcN(X)$ if and only if $x \nind^\mcN X$.
Also let $E$ be a $\es$-definable equivalence relation on $N$ with finitely many equivalence classes. 
Suppose that $Y$ is a minimal dependent set and $X$ an $E$-class such that $|Y| \geq 3$, 
$Y \setminus X \neq \es$, $Y \cap X \neq \es$ and $|Y \cap X|$ is minimal 
as $X$ ranges over the $E$-classes with which $Y$ has nonempty intersection.
Then there is a minimal dependent set $Z$ such that $|Z| \geq |Y|$, 
$|Z \cap X| = |Y \cap X| - 1$ and the number of $E$-classes with which $Z \setminus X$ has nonempty intersection 
is equal to the number of $E$-classes with which $Y \setminus X$ has nonempty intersection.
\end{lem}

\noindent
{\bf Proof.}
Suppose that $Y$ is a minimal dependent set of cardinality at least 3 which has nonempty intersection with at least two $E$-classes.
Let $X$ be an $E$-class such that $|Y \cap X|$ is minimal as $X$ ranges over the $E$-classes with which $Y$ has nonempty intersection.
Let $1 \leq k < l$ and let
\[
Y = \{y_1, \ldots, y_l\} \ \text{ where } \ Y \cap X = \{y_1, \ldots, y_k\}.
\]
Let $\bar{y}' = (y_1, \ldots, y_k)$, $y_i^0 = y_i$ for $i = 1, \ldots, l$ and $\bar{y}^0 = (y_{k+1}^0, \ldots, y_l^0)$.
By the existence of nonforking extensions (and $\omega$-saturation) there are, for $0 < n < \omega$,
$\bar{y}^n = (y_{k+1}^n, \ldots, y_l^n)$ such that for all $n < \omega$
\begin{equation*}
\bar{y}'\bar{y}^n \equiv_\mcN \bar{y}'\bar{y}^0 \ \text{ and } \
\bar{y}^{n+1} \underset{\bar{y}'}{\ind} \bar{y}^0 \ldots \bar{y}^n.
\end{equation*}
It follows that, for every $n < \omega$, $\bar{y}'\bar{y}^n$ is minimal dependent.
Since (by assumption) there are only finitely many $E$-classes, 
there are $s < t < \omega$ such that for every $k < i \leq l$,
$y_i^s$ and $y_i^t$ belong to the same $E$-class.
The lemma now follows if we can show that, for some proper subset
$Y^* \subset \{y_1, \ldots, y_k\}$, $Y^* \cup \rng(\bar{y}^s) \cup \rng(\bar{y}^t)$
is minimal dependent.
This follows from the following:

\bigskip

\noindent
{\bf Claim.} Suppose that $U$, $V$ and $W$ are mutually disjoint nonempty sets and that 
$U \cup V$ and $W \cup V$ are minimal dependent sets such that
$U \underset{V}{\ind} W$. Then there is a proper subset $V^* \subset V$ such that 
$U \cup V^* \cup W$ is minimal dependent.

\bigskip
\noindent
{\bf Proof of the claim.}
Let $v \in V$ and let $V' = V \setminus \{v\}$. 
Also let $u \in U$ and let $U' = U \setminus \{u\}$.
By assumption, $U \underset{V}{\ind} W$ so $U' \underset{V}{\ind} W$.
Since $UV$ is minimal dependent we have $U' \ind V$, so transitivity gives
$U' \ind VW$ and hence $U' \ind V'W$.
As $V' W$ is independent (since by assumption $VW$ is minimal dependent) it follows that
$U'V'W$ is independent. 
From the minimal dependence of $VW$ it follows that $v \in \cl_\mcN(V'W)$ and as $UV$ is minimal dependent
we get $u \in \cl_\mcN(U'V'W)$, so $UV'W$ is dependent and has dimension
$|U'V'W| = |UVW|-2$.

By arguing in the same way as above it follows that if $w \in W$ and $W' = W \setminus \{w\}$, 
then $UV'W'$ is independent and $UV'W \subseteq \cl_\mcN(UV'W')$.
If $V' = \es$ then $UW$ is minimal dependent and we are done, so now suppose that $V' \neq \es$.
If for any $v' \in V'$, and letting $V'' = V' \setminus \{v'\}$, the set $UV''W$ is independent then $UV'W$ is minimal dependent
and we are done.

Now suppose that $UV''W$ is dependent, where $V''$ is as above for some $v' \in V'$.
Then there is $S \subseteq UV''W$ such that $S$ is minimal dependent.
Since (as we have shown above) 
$U'V'W$ and $UV'W'$ are independent sets whenever $U' = U \setminus \{u\}$, $W' = W \setminus \{w\}$ for
$u \in U$ and $w \in W$ it follows that $UW \subseteq S$.
By letting $V^* = S \cap V''$ it follows that $UV^*W = S$ is minimal dependent.
This ends the proof of the claim and also of Lemma~\ref{reducing the number of types}.
\hfill $\square$
\\

\noindent
Now we can finish the proof of Theorem~\ref{nontrivial dependence implies nontrivial pregeometry}.
By Lemma~\ref{all elements have the same type}, there is an equivalence class $X$ of `$\sim$' such that
$X$ contains a set of cardinality at least 3 which is minimal dependent (over $\es$) with respect to $\mcM_P$.
By  Lemma~\ref{the structure on a single equivalence class}~(ii), this set is also minimal dependent with respect to $\mcM_P \uhrc X$.
This together with the other parts of the same lemma implies that if $\mcM' = \mcM_P \uhrc X$,
then conditions~(i) --~(iv) of 
Theorem~\ref{nontrivial dependence implies nontrivial pregeometry} hold and hence it is proved.

\section{Proof of Theorem~\ref{impossibility of a structure satisfying the other result}: 
impossibility of a nontrivial pregeometry (in the given context)}\label{impossibility of a nontrivial pregeometry}

\noindent
In order to prove Theorem~\ref{impossibility of a structure satisfying the other result} 
it suffices to derive a contradiction from the following assumptions which we now make.
Let $V$ be a ternary finite relational vocabulary and $\mcM$ a $V$-structure such that the following conditions hold:
\begin{itemize}
\item[(i)] $Th(\mcM)$ has elimination of quantifiers and is simple.

\item[(ii)] All elements of $M$ have the same type over $\acl_{(\mcM)\meq}(\es)$.

\item[(iii)] $(M, \cl)$ is a pregeometry if, for every $X \subseteq M$ and every $x \in M$, $x \in \cl(X)$ if and only if $x \nind X$, where $\nind$ is
with respect to $Th(\mcM)$.

\item[(iv)] $M$ has a minimal dependent subset of cardinality at least $3$.
\end{itemize}

\noindent
We first prove two lemmas and then derive a contradiction with the help of these.
For these lemmas the full assumption~(i) is not necessary, but it suffices to assume (besides ((ii)--(iv))
that $Th(\mcM)$ is simple and that $\mcM$ is $\omega$-saturated.

\begin{lem}\label{every set is contained in a minimal dependent set}
If $A \subseteq M$ is finite and independent then there is a minimal dependent 
$A' \subseteq M$ such that $A \subset A'$ and $|A'| \geq |A|+2$.
\end{lem}

\noindent
{\bf Proof.}
We prove the lemma by induction on $n = |A|$.
If $n = 1$ then the conclusion follows from assumptions~(ii) and~(iv).
So now suppose that $n \geq 1$ and $\{a_1, \ldots, a_n, b_1\} \subseteq M$ is independent.
By the induction hypothesis there is $a_{n+1} \in M$ such that $\{a_1, \ldots, a_n, a_{n+1}\}$ is independent and included in
a minimal dependent set.
By the existence of nonforking extensions we may assume that $a_{n+1} \underset{a_1 \ldots a_n}{\ind} b_1$.
Since $\{a_1, \ldots, a_n, b_1\}$ is independent it follows from transitivity that $\{a_1, \ldots, a_{n+1}, b_1\}$ is independent.
Choose $b_{n+1} \in M$ such that $a_1a_{n+1} \equiv_\mcM b_1b_{n+1}$.
Note that we now have $a_1\ldots a_n \ind b_1$, $a_1\ldots a_n \ind a_{n+1}$ and $b_1 \ind b_{n+1}$.
By assumption~(ii), all elements have the same type over $\acl_{\mcM\meq}(\es)$, so 
(as $\mcM$ is $\omega$-categorical) the independence theorem of simple theories implies
that there is $c_{n+1} \in M$ such that 
\begin{align}\label{c-n+1 has the right types}
&a_1\ldots a_n a_{n+1} \equiv_\mcM a_1\ldots a_n c_{n+1}, \\
&b_1 b_{n+1} \equiv_\mcM b_1 c_{n+1} \ \ \text{ and} \nonumber \\
&c_{n+1} \ind a_1\ldots a_n b_1. \nonumber
\end{align}
As $\{a_1, \ldots, a_n, b_1\}$ is independent (by assumption) it follows that 
\begin{equation}\label{a-1 to a-n, b-1 and c-n+1}
\{a_1, \ldots, a_n, c_{n+1}, b_1\} \text{ is independent, so in particular } 
a_1 \ldots a_n c_{n+1} \ind b_1.
\end{equation}
By the choice of $a_{n+1}$ and~(\ref{c-n+1 has the right types}),
there are $m \geq n+2$ and $a_{n+2}, \ldots, a_m \in M$ such that 
\begin{equation*}
\{a_1, \ldots, a_n, c_{n+1}, a_{n+2}, \ldots, a_m\} \text{ is minimal dependent.}
\end{equation*}
By the existence of nonforking extensions we may assume that
\begin{equation*}
a_{n+2} \ldots a_m \underset{a_1 \ldots a_n c_{n+1}}{\ind} b_1.
\end{equation*}
This together with~(\ref{a-1 to a-n, b-1 and c-n+1}) and transitivity implies that
\begin{equation}\label{b-1 is independent from c-n+1 and the a}
a_1 \ldots a_n c_{n+1} a_{n+2} \ldots a_m \ind b_1.
\end{equation}
As $b_1 c_{n+1} \equiv_\mcM b_1 b_{n+1} \equiv_\mcM a_1 a_{n+1}$ there are $b_2, \ldots,$ $b_n \in M$
and $b_{n+2}, \ldots, b_m \in M$ such that 
\begin{equation*}
\{b_1, \ldots, b_n, c_{n+1}, b_{n+2}, \ldots, b_m\} \text{ is minimal dependent.}
\end{equation*}
By the existence of nonforking extensions we may assume that 
\begin{equation}\label{independence of all b from all a over}
b_2 \ldots b_n b_{n+2} \ldots b_m \underset{b_1 c_{n+1}}{\ind} a_1 \ldots a_n a_{n+2} \ldots a_m.
\end{equation}
Let 
\begin{equation*}
A = \{a_1, \ldots, a_n, a_{n+2}, \ldots, a_m\} \ \text{ and } \  B = \{b_1, \ldots, b_n, b_{n+2}, \ldots, b_m\}.
\end{equation*}

It now suffices to prove that $AB$ is minimal dependent, because $a_1,$ $\ldots,$ $a_n, b_1 \in AB$ and
$|AB| = 2m-2 \geq |\{a_1, \ldots, a_n, b_1\}| + 2$, since $m \geq n+2$.
Since $Ac_{n+1}$  and $Bc_{n+1}$ are minimal dependent it follows that $c_{n+1} \in \cl(A)$ 
and $b_m \in \cl(Bc_{n+1} \setminus \{b_m\})$.
Consequently (as $(M, \cl)$ is a pregeometry) $b_m \in \cl(AB \setminus \{b_m\})$, so $AB$ is dependent.
It remains to prove that every proper subset of $AB$ is independent.
We will divide the argument into a few cases.

Let $A' \subset A$ be a proper subset (and we will show that $A'B$ is independent).
Then $A' \ind c_{n+1}$.
From~(\ref{b-1 is independent from c-n+1 and the a}) we get (by monotonicity)
$A' \underset{c_{n+1}}{\ind} b_1$,
so transitivity gives $A' \ind b_1 c_{n+1}$.
This together with~(\ref{independence of all b from all a over}) and transitivity gives
\[
A' \ind b_1 \ldots b_n c_{n+1} b_{n+2} \ldots b_n
\]
and hence $A' \ind B$ from which it follows (by repeatedly using monotonicity and transitivity) that
$A'B$ is independent.
Since we can choose $A' \subset A$ so that $|A'| = |A|-1$ it follows that $AB$ has dimension 
$|AB|-1$ (where dimension is with respect to `$\cl$').

Now let $B' \subset B$ with $|B'| = |B|-1$.
To complete the proof it suffices to prove that $AB'$ is independent.
There is $b \in B$ such that $B' = B \setminus \{b\}$.
First suppose that $b \neq b_1$, so $b_1 \in B'$.
As $Bc_{n+1}$ is minimal dependent, $B'c_{n+1}$ is independent,
so $b_1c_{n+1} \ind (B' \setminus \{b_1\})$.
This together with~(\ref{independence of all b from all a over}) and transitivity (and monotonicity) implies that
\[
(B' \setminus \{b_1\}) \ind b_1 a_1 \ldots a_n a_{n+2} \ldots a_m
\]
and hence
\[
(B' \setminus \{b_1\}) \underset{b_1}{\ind} a_1 \ldots a_n a_{n+2} \ldots a_m.
\]
This together with~(\ref{b-1 is independent from c-n+1 and the a}) and transitivity gives
$B' \ind A$, so $AB'$ is independent.

Now suppose that $b = b_1$, so $B' = B \setminus \{b_1\}$.
Towards a contradiction suppose that $AB'$ is dependent.
Then there is a proper subset $C \subset AB'$ such that
$AB' \subseteq \cl(C)$.
As $Ac_{n+1}$ is minimal dependent we have $c_{n+1} \in \cl(A)$.
As $Bc_{n+1}$ is minimal dependent we have $b_1 \in \cl(B'c_{n+1})$ and hence
$b_1 \in \cl(AB')$. Since $AB' \subseteq \cl(C)$ we get $AB = AB'b_1 \subseteq \cl(C)$
where $|C| \leq |AB| - 2$. 
Thus the dimension of $AB$ is at most $|AB| - 2$. 
But this contradicts our earlier conclusion that the dimension of $AB$ is $|AB| - 1$.
\hfill $\square$

\begin{lem}\label{turning a minimal dependent set into an independent set}
Let $A \subset M$ be a minimal dependent set of cardinality at least $3$ and let $a_1, a_2, a_3 \in A$ be distinct.
Then there is $a'_3 \in M$ such that $(A \setminus \{a_3\}) \cup \{a'_3\}$ is independent and
for all $b, c \in A \setminus \{a_3\}$, if $\{b, c\} \neq \{a_1, a_2\}$, then $bca'_3 \equiv_\mcM bca_3$.
\end{lem}

\noindent
{\bf Proof.}
Let $A = \{a_1, \ldots, a_n\} \subset M$ be a minimal dependent where $n \geq 3$.
Let $B = \{a_4, \ldots, a_n\}$ (or $B = \es$ if $n = 3$).
Then $a_1 \underset{B}{\ind} a_2$, $a_3 \underset{B}{\ind} a_1$ and $a_3 \underset{B}{\ind} a_2$.
By the independence theorem of simple theories, there is a type $q$ over 
$\acl_{\mcM\meq}(B) \cup \{a_1, a_2\}$ which extends the type of $a_3$ over 
$\acl_{\mcM\meq}(B) \cup \{a_1\}$ and the type of $a_3$ over $\acl_{\mcM\meq}(B) \cup \{a_2\}$, and $q$ does not 
divide over $B$.
Since $Th(\mcM)$ is $\omega$-categorical, and hence $\mcM\meq$ is $\omega$-saturated, there is $a'_3 \in M$ which realizes $q$, so 
$\tp_\mcM(a'_3 / a_1 B) = \tp\mcM(a_3 / a_1 B)$, $\tp_\mcM(a'_3 / a_2 B) = \tp_\mcM(a_3 / a_2 B)$ and
$a'_3 \underset{B}{\ind} a_1 a_2$.
It follows, since $A$ is minimal dependent, that $\{a_1, a_2, a'_3, a_4, \ldots, a_n\} = \{a'_3, a_1, a_2\} \cup B$ is independent.
From the choice of $a'_3$  it also follows that if 
$b, c \in A \setminus \{a_3\}$ and $\{b, c\} \neq \{a_1, a_2\}$, then $bca'_3 \equiv_\mcM bca_3$.
\hfill $\square$
\\

\noindent
Now we are ready to make a construction which will lead to a contradiction.
By the existence of nonforking extensions we find $a_1, a_2, a_3 \in M$ such 
that $A_0 = \{a_1, a_2, a_3\}$ is independent.
Let $a_3^0 = a_3$.
We now find elements $a_3^n \in M$, for $0 < n < \omega$, and finite sets
$A_n, B_n \subseteq M$, for $n < \omega$, such that for every $n < \omega$ the following conditions hold:
\begin{itemize}
\item[(a)] $A_n = \{a_1, a_2, a_3^n\}$, $A_{n+1} = \{a_1, a_2, a_3^{n+1}\}$ and both $A_n$ and $A_{n+1}$ are independent,

\item[(b)] $A_n B_0 \ldots B_n$ is minimal dependent,

\item[(c)] $A_{n+1}B_0 \ldots B_n$ is independent, and

\item[(d)] if $b, c \in A_{n+1}B_0 \ldots B_n \setminus \{a_3^{n+1}\}$ and $\{b, c\} \neq \{a_1, a_2\}$, 
then $bca_3^{n+1} \equiv_\mcM bca_3^n$.
\end{itemize}
Suppose that (a)--(d) hold.
By Lemma~\ref{every set is contained in a minimal dependent set}
there is $B_{n+1}$ such that $A = A_{n+1}B_0 \ldots B_{n+1}$ is minimal dependent.
By~(a) and
Lemma~\ref{turning a minimal dependent set into an independent set}
there is $a_3^{n+2} \in M$ such that
\begin{itemize}
\item $(A_{n+1}B_0 \ldots B_{n+1} \setminus \{a_3^{n+1}\}) \cup \{a_3^{n+2}\}$ is independent and
\item if $b, c \in A_{n+1}B_0 \ldots B_{n+1} \setminus \{a_3^{n+1}\}$ and $\{b, c\} \neq \{a_1, a_2\}$, 
then \\
$bca_3^{n+2} \equiv_\mcM bca_3^{n+1}$.
\end{itemize}
Let $A_{n+2} = \{a_1, a_2, a_3^{n+2}\}$.
Now~(a)--(d) holds when `$n$' is replaced by `$n+1$'.
Note that if $n=0$ then~(a) holds by the choices of $A_0$ and $a_3^0$ and then we find $B_0$, $A_1$ and $a_3^1$ such that~(b)--(d) hold
for $n=0$ in the same way as we did for the general case $n$.

Since $Th(\mcM)$ is $\omega$-categorical (because $\mcM$ has elimination of quantifiers)
there are only finitely many 3-types over $\es$.
Thus there are $i < j$ such that $a_1a_2a_3^{i+1} \equiv_\mcM a_1a_2a_3^{j+1}$.
Since the vocabulary is ternary it follows from elimination of quantifiers and (d) applied to
$n = i+1, \ldots, j$ that
\begin{equation}\label{same type for i and j}
\tp_\mcM(a_1, a_2, a_3^{i+1} / B_0 \ldots B_{i+1}) = \tp_\mcM(a_1, a_2, a_3^{j+1} / B_0 \ldots B_{i+1}).
\end{equation}
By~(b) for $n = i$, $A_{i+1}B_0 \ldots B_{i+1}$ is dependent.
By~(c) for $n = j$, \\
$A_{j+1}B_0 \ldots B_j$ is independent and as $i < j$ it follows that $A_{j+1}B_0 \ldots B_{i+1}$ is independent.
Since (by~(a) for $n = i+1$ and $n = j+1$) $A_{i+1} = \{a_1, a_2, a_3^{i+1}\}$ and $A_{j+1} = \{a_1, a_2, a_3^{j+1}\}$ we must have
\[
\tp_\mcM(a_1, a_2, a_3^{i+1} / B_0 \ldots B_{i+1}) \neq \tp_\mcM(a_1, a_2, a_3^{j+1} / B_0 \ldots B_{i+1})
\]
which contradicts~(\ref{same type for i and j}).
Thus the proof of Theorem~\ref{impossibility of a structure satisfying the other result} 
is finished.

\begin{rem}\label{remark about higher dimensional amalgamation}{\rm
The notion of ``$n$-amalgamation property'' has been defined in slightly different ways and with slightly different
names in various articles.
Let us use the definition of {\em $n$-complete amalgamation (property)} used in
\cite{KKT, Pal17}. Then the independence theorem is equivalent to the $3$-complete amalgamation property.
Now one can modify Lemma~\ref{turning a minimal dependent set into an independent set}, its proof and
the argument after it so that one gets the following result:
{\em Suppose that $V$ is a finite relational vocabulary with maximal arity $k$ and $T$ a simple $V$-theory
with elimination of quantifiers and $k$-complete amalgamation. 
Then $T$ has no model $\mcM$ such that (ii) -- (iv) in the beginning of this section hold.}
Since $n$-complete amalgamation, when restricted to (finite tuples of) real elements, is preserved when passing from $T$ to 
$Th(\mcM')$ in Theorem~\ref{nontrivial dependence implies nontrivial pregeometry} it follows that if
$T$ is as above then $T$ has trivial dependence and finite SU-rank.
}\end{rem}

\noindent
{\bf Acknowledgement.} I thank the anonymous referee for a detailed examination of the article,
including finding some minor mistakes (now corrected) and
giving suggestions that improved the clarity of the arguments.


\begin{thebibliography}{99}\label{References}

\bibitem{Cas} E. Casanovas, {\em Simple Theories and Hyperimaginaries}, Cambridge University Press (2011).

\bibitem{CH} G. Cherlin, E. Hrushovski, Finite Structures with Few Types, 
{\em Annals of Mathematics Studies}, Nr. 152, Princeton University Press (2003).

\bibitem{CHL} G. Cherlin, L. Harrington, A. H. Lachlan, $\omega$-categorical $\omega$-stable structures,
{\em Annals of Pure and Applied Mathematics}, Vol. 28 (1986) 103--135.

\bibitem{DK} T. de Piro, B. Kim, The geometry of 1-based minimal types,
{\em Transactions of the American Mathematical Society}, Vol. 355 (2003) 4241--4263.

\bibitem{EW} D. Evans, F. O. Wagner, Supersimple $\omega$-categorical groups and theories,
{\em The Journal of Symbolic Logic}, Vol. 65 (2000) 767--776.

\bibitem{Goode} J. B. Goode, Some trivial considerations, {\em The Journal of Symbolic Logic},
Vol. 56 (1991) 624--631.

\bibitem{HKP} B. Hart, B. Kim, A. Pillay, Coordinatisation and canonical bases in simple theories,
{\em The Journal of Symbolic Logic}, Vol. 65 (2000) 293--309.

\bibitem{Hod} W. Hodges,  {\em Model theory}, Cambridge University Press (1993).

\bibitem{Hru85} E. Hrushovski, Locally modular regular types,
in J. T. Baldwin (ed.) {\em Classification Theory, Proceedings, Chicago 1985}, Springer-Verlag, Berlin.

\bibitem{Hru_pseudoplane} E. Hrushovski, A stable $\aleph_0$-categorical pseudoplane, unpublished notes (1988).

\bibitem{Hru90} E. Hrushovski, Unidimensional theories are superstable,
{\em Annals of Pure and Applied Logic}, Vol. 50 (1990) 117--138.

\bibitem{Hru92} E. Hrushovski, Unimodular minimal theories,
{\em Journal of The London Mathematical Society, Ser. 2} Vol. 46 (1992) 385--396.

\bibitem{Kim98} B. Kim, A note on Lascar strong types in simple theories,
{\em The Journal of Symbolic Logic}, Vol. 63 (1998) 926--936.

\bibitem{Kim_book} B. Kim, {\em Simplicity Theory}, Oxford University Press (2014).

\bibitem{KKT} B. Kim, A. S. Kolesnikov, A. Tsuboi, Generalized amalgamation and $n$-simplicity,
{\em Annals of Pure and Applied Logic}, Vol. 155 (2008) 97--114.

\bibitem{Kop16PAM} V. Koponen, Binary simple homogeneous structures are supersimple with finite rank,
{\em Proceedings of the American Mathematical Society}, Vol. 144 (2016) 1745--1759.

\bibitem{KopBin} V. Koponen, Binary simple homogeneous structures, 
{\em Annals of Pure and Applied Logic}, to appear, online\\
\url{https://arxiv.org/abs/1609.02433}

\bibitem{Kop17a} V. Koponen, On constraints and dividing in ternary homogeneous structures, (submitted),
online:\\
\url{https://arxiv.org/abs/1707.05954}

\bibitem{Lach74} A. H. Lachlan, Two conjectures regarding the stability of $\omega$-categorical theories,
{\em Fundamenta Mathematicae}, Vol. 81 (1974) 133 -- 145.

\bibitem{Lach97} A. H. Lachlan, Stable finitely homogeneous structures: a survey,
in B. T. Hart et. al. (eds.), {\em Algebraic Model Theory}, 145--159, Kluwer Academic Publishers (1997) 

\bibitem{Mac91} D. Macpherson, Interpreting groups in $\omega$-categorical structures,
{\em The Journal of Symbolic Logic}, Vol. 56 (1991) 1317--1324.

\bibitem{Mac11} D. Macpherson, A survey of homogeneous structures,
{\em Discrete Mathematics}, Vol. 311 (2011) 1599--1634.

\bibitem{Pal17} D. Palac\'{i}n, {\em Generalized amalgamation and homogeneity}, 
{\em The Journal of Symbolic Logic}, Vol. 82 (2017) 1402--1421.

\bibitem{TZ} K. Tent, M. Ziegler, {\em A course in model theory}, Lecture Notes in Logic 40,
Cambridge University Press (2012).

\bibitem{Wag} F. O, Wagner, {\em Simple Theories}, Kluwer Academic Publishers (2000).

\bibitem{Zilber_mono} B. Zilber, {\em Uncountably Categorical Theories}, AMS Translations of Mathematical Monographs,
Vol. 117 (1993).

\end{thebibliography}
\end{document}